\documentclass[preprint,12pt]{elsarticle}

\usepackage{amssymb}
\usepackage{amsmath,amssymb,amsfonts}
\usepackage{tikz}

\usetikzlibrary{decorations.pathmorphing, patterns, decorations.markings}

\newdefinition{rmk}{Remark}
\newproof{proof}{Proof}
\journal{ArXiv}

\begin{document}


\begin{frontmatter}

\title{The continuity properties of the solutions of the wave equation at different boundary conditions across the characteristic line}

\author[1]{Alemdar Hasanov\corref{cor1}
\fnref{fn1}}
\ead{alemdar.hasanoglu@gmail.com}
\author[2]{Vladimir G. Romanov \fnref{fn2}}
\ead{romanov0511@gmail.com}
\cortext[cor1]{Corresponding author}
\fntext[fn1]{Emeritus Professor, Kocaeli University,  T\"{u}rkiye}
\fntext[fn2] {Sobolev Institute of Mathematics, Russia}
\address[1]{Kocaeli University, T\"{u}rkiye}
\address[2]{Sobolev Institute of Mathematics, Russia}

\begin{abstract}
It is well-known the continuity properties of solutions to the wave equation depend strongly on both the initial data and the boundary conditions, particularly across the characteristic lines, where information propagates. Moreover, the continuity properties of solutions to the wave equation depend strongly on the type of boundary condition imposed and on the interaction of the characteristics with the boundary. Across the characteristic lines, singularities and discontinuities propagate according to the geometry of the problem, while the boundary conditions determine how these singularities are reflected or transmitted.

In this study, we will examine the behavior of solutions of the wave equation $u_{tt}=c^2 u_{xx}$, $c>0$, corresponding to different boundary conditions, on its characteristic line, using the formulas obtained as a result of d'Alembert's formula, for all the three, Dirichlet, Neumann and Robin, types of boundary conditions.

\end{abstract}
\begin{keyword}
Wave equation, d'Alembert's formuls, characteristics, continuity across the characteristics.
\end{keyword}

\end{frontmatter}


\section{Introduction}

The behavior of solutions to boundary value problems (BVPs) related to the wave equation across the characteristic line is outlined in \cite{Courant:1953}, \cite{Evans:2010} and \cite{FJohn:1978}. However, the inverse coefficient problems for wave equation formulated and examined in detail in articles \cite{VGR-AH-1:2020}, \cite{VGR-AH-2:2020} and \cite{VGR-AH-3:2021} necessitate a more detailed investigation of the continuity property of the solution of the corresponding forward problem across the characteristic line. Based on this demand, the investigation of the continuity of solutions corresponding to different boundary conditions (Dirichlet, Neumann and Robin) across the characteristics is addressed in this study. Namely, we will examine the behavior of solutions of the wave equation $u_{tt}=c^2 u_{xx}$, $c>0$, corresponding to different boundary conditions, on its characteristic line, using the formulas obtained as a result of d'Alembert's formula, for all the three, Dirichlet, Neumann and Robin, types of boundary conditions.

\section{d'Alembert's solution and applications}

Consider the problem
\begin{eqnarray}\label{1-1}
\left\{\begin{array}{ll}
	u_{tt}-c^2 u_{xx}=0,~ (x,t)\in \mathbb{R} \times (0,+\infty),\\ [6pt]
    u(x,0)=g(x), ~ u_t(x,0)=h(x), ~ x\in \mathbb{R},
\end{array}\right.
\end{eqnarray}
where $c>0$ and $\mathbb{R} :=(-\infty,+\infty)$.
Assume that
\begin{eqnarray}\label{1-2}
g\in C^2(\mathbb{R}), ~ h\in C^1(\mathbb{R}).
\end{eqnarray}
Then the solution $u(x,t)$ of the problem (\ref{1-1}) belongs to $C^2(\mathbb{R}\times (0,+\infty))$, and is defined by d'Alembert's formula \cite{Evans:2010}:
\begin{eqnarray}\label{1-3}
u(x,t)=\frac{1}{2}\,\left [ g(x+ct)+g(x-ct)\right ]	+\frac{1}{2c}\,\int_{x-ct}^{x+ct} h(\xi)d\xi,\, (x,t)\in \mathbb{R} \times (0,+\infty).
\end{eqnarray}

\subsection{The Dirichlet boundary condition}

To illustrate the first application of d'Alembert's formula we consider the following initial boundary value problem (IBVP):
\begin{eqnarray}\label{1-4}
\left\{\begin{array}{ll}
	u_{tt}-c^2 u_{xx}=0,~ (x,t)\in \mathbb{R}_+ \times (0,+\infty),\\ [6pt]
    u(x,0)=g(x), ~ u_t(x,0)=h(x), ~ x\in \mathbb{R}_+,\\ [6pt]
    u(0,t)=\nu(t), \,t \in (0,+\infty),
\end{array}\right.
\end{eqnarray}
where $\mathbb{R}_+ :=(0,+\infty)$.

It is assumed that the Dirichlet boundary data satisfy the regularity condition $\nu \in C^2(0,+\infty)$  to ensure the existence of a classical solution to problem (\ref{1-4}).

  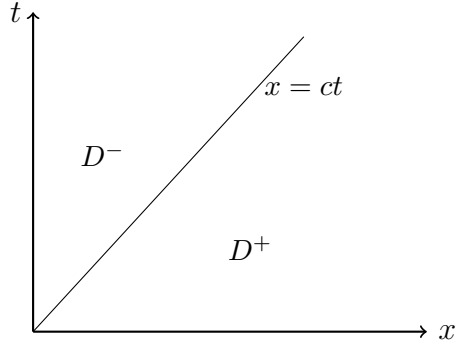
\begin{figure}
  $\quad$\\
  $\quad$\\
  \centering
 \hspace*{-.6cm} \begin{tikzpicture}[scale=.65]
  \draw[thick,->] (0,0) -- (8.0,0) node [right]{$x$};
      \draw[thick,->] (0,0) -- (0,6.5) node [left]{$t$};
      \draw[color=black] (0,0) -- (5.5,6.0);
       \node[label=left:{\small$D^{+}$}] at (5.3,1.7) {};
       \node[label=left:{\small$D^{-}$}] at (2.3,3.6) {};
       \coordinate (P) at (5.5,5) ;
\draw[rotate around={45:(P)}] (P) node {\small$x=ct$} ;
        \end{tikzpicture}
  \caption{Domains $D^{-}$ and $D^{+}$ formed by the characteristic line $x-ct=0$}
  \label{Fig-1}
  \end{figure}

We introduce the domains
\begin{eqnarray}\label{1-5}
\left. \begin{array}{ll}
D^{+}:=\left \{(x,t)\in \mathbb{R}_+\,: \, x > ct,\, t>0 \right \}, \\ [6pt]
D^{-}:=\left \{(x,t)\in \mathbb{R}_+\,: \, 0<x < ct,\, t>0 \right \},
\end{array}\right.
\end{eqnarray}
formed by characteristic line $x-ct=0$ (Figure \ref{Fig-1}).

Evidently, in the domain $D^{+}$ the solution $u(x,t)$ of problem (\ref{1-4}) is defined uniquely by d'Alembert's formula (\ref{1-3}). What about the domain $D^{-}\,$? Defining the Dirichlet boundary condition, in a sense, necessary to find the answer to this question.

Indeed, the general solution of problem (\ref{1-1}) is
\begin{eqnarray}\label{1-6}
u(x,t)=\phi(x+ct)+\psi(x-ct),
\end{eqnarray}
and the initial conditions in (\ref{1-4}) imply that
\begin{eqnarray*}
\left\{\begin{array}{ll}
	\phi(x)+\psi(x)=g(x),\\ [6pt]
\displaystyle \phi(x)-\psi(x)=\frac{1}{c}\int_0^x h(\xi)d\xi.
\end{array}\right.
\end{eqnarray*}
Hence,
\begin{eqnarray}\label{1-7}
\left\{\begin{array}{ll}
\displaystyle	\phi(x)=\frac{1}{2}\,g(x)+ \frac{1}{2c}\int_0^x h(\xi)d\xi,\\ [14pt]
\displaystyle \psi(x)=\frac{1}{2}\,g(x)- \frac{1}{2c}\int_0^x h(\xi)d\xi.
\end{array}\right.
\end{eqnarray}
Therefore, $\phi(x+ct)$ and $\psi(x-ct)$  are defined respectively as $x+ct > 0$
(which holds automatically as $x > 0,\, t > 0$), and as $x-ct > 0$ (which is holds only
for $x > ct$), that is in $D^{+}$. To define $\psi(x-ct)$  for $0<x<ct$, i.e. in the domain $D^{-}$, we need to define $\psi(y)$ for $y<0$.

To this end, we employ the boundary condition $u(0,t)=\nu(t)$ in (\ref{1-6}). This yields:
\begin{eqnarray*}
\phi(ct)+\psi(-ct)=\nu(t), t>0~,
\end{eqnarray*}
or equivalently, $\phi(-x)+\psi(x)=\nu(-x/c)$ for $x<0$ (and $-x>0$), where we used $t:=-x/c$. Hence $\psi(x)=\nu(-x/c)-\phi(-x)$, and with the first formula of (\ref{1-7}) this implies:
\begin{eqnarray}\label{1-8}
\psi(x)=\nu(-x/c)-\phi(-x) \qquad  \qquad \qquad \qquad \qquad \qquad \qquad \quad \nonumber \\  [2pt]
\qquad \qquad = \nu(-x/c)-\frac{1}{2}\,g(-x)- \frac{1}{2c}\int_0^{-x} h(\xi)d\xi,~x<0,\, t>0.
\end{eqnarray}
Replace $x$ with $x+ct$ in the first formula of (\ref{1-7}) and $x$ with $x-ct$ in formula (\ref{1-8}). With d'Alembert's formula (\ref{1-3}) this implies:
\begin{eqnarray}\label{1-9}
\displaystyle u(x,t)=\frac{1}{2}\,g(x+ct)+\frac{1}{2c}\int_0^{x+ct} h(\xi)d\xi\qquad \qquad \qquad \qquad \qquad \qquad \qquad \nonumber\\ [2pt]
\qquad \qquad -\frac{1}{2}\,g(-x+ct)-\frac{1}{2c}\,\int_0^{ct-x} h(\xi)d\xi+\nu(t-x/c),~ (x,t)\in D^{-}.
\end{eqnarray}
Recall that the solution of the IBVP (\ref{1-4}) in the domain $D^{+}$ is given by formula (\ref{1-3}).

Thus, the solution of the IBVP (\ref{1-4}) obtained from d'Alembert's formula (\ref{1-3}) is
\begin{eqnarray}\label{1-10}
u(x,t)=
\left\{\begin{array}{ll}
\displaystyle \frac{1}{2}\,\left [ g(x+ct)+g(x-ct)\right ]	+\frac{1}{2c}\,\int_{x-ct}^{x+ct} h(\xi)d\xi,\,(x,t)\in D^{+},\\ [16pt]
\displaystyle \frac{1}{2}\,\left [ g(x+ct)-g(-x+ct)\right ]+\frac{1}{2c}\,\left [\int_0^{x+ct} h(\xi)d\xi -\int_0^{ct-x} h(\xi)d\xi\right ]\\ [16pt]
\qquad \qquad \qquad \qquad \qquad \qquad \qquad \qquad \qquad +\nu(t-x/c),\,(x,t)\in D^{-}\,.
\end{array}\right.
\end{eqnarray}

If $\nu(t)=0$ and $c=1$, this formıla is the same as the formula (10) in \cite{Evans:2010} (Ch. 2.4), obtained through the reflection method.

\subsection{The Neumann boundary condition}

As a second application of d'Alembert's formula consider the following IBVP:
\begin{eqnarray}\label{1-11}
\left\{\begin{array}{ll}
	u_{tt}-c^2 u_{xx}=0,~ (x,t)\in \mathbb{R}_+ \times (0,+\infty),\\ [6pt]
    u(x,0)=g(x), ~ u_t(x,0)=h(x), ~ x\in \mathbb{R}_+,\\ [6pt]
    u_x(0,t)=f(t), \,t \in (0,+\infty).
\end{array}\right.
\end{eqnarray}

It is assumed that the Neumann boundary data satisfy the regularity condition $f \in C^{1}(0,+\infty)$ to ensure the existence of a classical solution to problem (\ref{1-11}).

In view of formula (\ref{1-6}) the Neumann condition $u_x(0,t)=f(t)$ yields: $\phi'(ct)+\psi'(-ct)=f(t)$, as $t>0$, or equivalently,
\begin{eqnarray*}
\phi(-x)-\psi(x)=c\,\int_0^{-x/c} f(\xi)d\xi,
\end{eqnarray*}
for $x<0$ (and $-x>0$), where we used $t:=-x/c$. With the first formula of (\ref{1-7}) this implies:
\begin{eqnarray}\label{1-12}
\psi(x)=-c\,\int_0^{-x/c} f(\xi)d\xi+\phi(-x) \qquad  \qquad \qquad \qquad \qquad \qquad  \qquad \nonumber \\  [2pt]
\qquad \qquad  \qquad = -c\,\int_0^{-x/c} f(\xi)d\xi+\frac{1}{2}\,g(-x)+\frac{1}{2c}\,\int_0^{-x} h(\xi)d\xi,~t>0.
\end{eqnarray}
Again, replacing $x$ with $x+ct$ in the first formula of (\ref{1-7}) and $x$ with $x-ct$ in formula (\ref{1-12}), and taking formula (\ref{1-3}) into account we conclude that
\begin{eqnarray*}
\displaystyle u(x,t)=\frac{1}{2}\,g(x+ct)+\frac{1}{2c}\int_0^{x+ct} h(\xi)d\xi\qquad \qquad \qquad \qquad \qquad  \qquad \qquad \qquad \quad \nonumber\\ [4pt]
\qquad \qquad +\frac{1}{2}\,g(-x+ct)+\frac{1}{2c}\,\int_0^{-x+ct} h(\xi)d\xi-c \,\int_0^{t-x/c} f(\xi)d\xi ,~ (x,t)\in D^{-}.
\end{eqnarray*}
Hence, the solution of the IBVP (\ref{1-11}) obtained from d'Alembert's formula (\ref{1-3}) is
\begin{eqnarray}\label{1-13}
u(x,t)=
\left\{\begin{array}{ll}
\displaystyle \frac{1}{2}\,\left [ g(x+ct)+g(x-ct)\right ]	+\frac{1}{2c}\,\int_{x-ct}^{x+ct} h(\xi)d\xi,\,(x,t)\in D^{+},\\ [16pt]
\displaystyle \frac{1}{2}\,\left [ g(x+ct)+g(-x+ct)\right ]+\frac{1}{2c}\,\left [\int_0^{x+ct} h(\xi)d\xi +\int_0^{ct-x} h(\xi)d\xi\right ] \\ [16pt]
\displaystyle \qquad \qquad \qquad \qquad \qquad  \qquad \qquad \qquad
-c \int_0^{t-x/c} f(\xi)d\xi,\, (x,t)\in D^{-}\,.
\end{array}\right.
\end{eqnarray}

It can be easily verified that the solution $u(x,t)$ defined by (\ref{1-13}) satisfies the equation, initial and boundary conditions in problem (\ref{1-11}).

\subsection{The Robin boundary condition}

Finally, consider the following IBVP:
\begin{eqnarray}\label{1-14}
\left\{\begin{array}{ll}
	u_{tt}-c^2 u_{xx}=0,~ (x,t)\in \mathbb{R}_+ \times (0,+\infty),\\ [6pt]
    u(x,0)=g(x), ~ u_t(x,0)=h(x), ~ x\in \mathbb{R}_+\,,\\ [6pt]
    u_x(0,t)+\sigma u(0,t)=f(t), \,t \in (0,+\infty),
\end{array}\right.
\end{eqnarray}
where $\sigma>0$.

In view of the Robin condition in (\ref{1-14}) and formula (\ref{1-6}) we deduce that
\begin{eqnarray}\label{1-15}
\psi'(-ct) +\sigma \psi(-ct) +\phi'(ct) +\sigma \phi(ct)=f(t), \, t>0.
\end{eqnarray}
With the initial condition $\psi(0)=g(0)/2$, this leads to the Cauchy problem
\begin{eqnarray}\label{1-16}
\left\{\begin{array}{ll}
\psi'(x) +\sigma \psi(x)=F(-x), \, t>0.\\ [6pt]
\displaystyle \psi(0)= \frac{1}{2}\,g(0),
\end{array}\right.
\end{eqnarray}
with
\begin{eqnarray}\label{1-17}
F(-x)=f(-x/c)-\phi'(-x)-\sigma \phi(-x), \, t>0.
\end{eqnarray}
The solution of this problem is
\begin{eqnarray}\label{1-18}
\displaystyle \psi(x)= \frac{1}{2}\,g(0)\exp(-\sigma x)+\int_0^x \exp(-\sigma (x-\xi))\,F(-\xi)\,d\xi, \, x<0.
\end{eqnarray}
Substituting (\ref{1-17}) in (\ref{1-18}) we obtain:
\begin{eqnarray*}
\displaystyle \psi(x)= \frac{1}{2}\,g(0)\exp(-\sigma x)+\int_0^x \exp(-\sigma (x-\xi))\,f(-\xi/c)\,d\xi \qquad \qquad \qquad \qquad \nonumber \\ [6pt]
\qquad \qquad -2\sigma \int_0^x \exp(-\sigma (x-\xi))\,\phi(-\xi)\,d\xi+\phi(-\xi)-\phi(0)\exp(-\sigma x), \, x<0.
\end{eqnarray*}

Taking into account formula (\ref{1-6}), and then employing the above methodology we get:
\begin{eqnarray*}
\displaystyle u(x,t)= \frac{1}{2}\,g(x+ct)+\frac{1}{2c}\,\int_0^{x+ct} h(\xi)d\xi \qquad \qquad \qquad \qquad \qquad \qquad \qquad \qquad \qquad \qquad\nonumber \\ [6pt]
+\frac{1}{2}\,g(0)\exp(-\sigma (x-ct))+\exp(-\sigma (x-ct)) \int_0^{x-ct} \exp(\sigma \xi)\,f(-\xi/c)\,d\xi \qquad \qquad \qquad \nonumber \\ [6pt]
-2\sigma \exp(-\sigma (x-ct)) \int_0^{x-ct} \exp(\sigma \xi)\,\phi(-\xi)\,d\xi+\frac{1}{2}\,g(-x+ct)+\frac{1}{2c}\,\int_0^{ct-x} h(\xi)d\xi \quad  \nonumber \\ [6pt]
+\frac{1}{2}\,g(0)\exp(-\sigma (x-ct)), ~ (x,t)\in D^{-}\,.
\end{eqnarray*}
Eliminating identical terms and grouping them, we finally arrive at the following formula for the solution of the Robin problem (\ref{1-14}) in the domain $D^{-}$:
\begin{eqnarray*}
\displaystyle u(x,t)= \frac{1}{2}\,\left [g(x+ct)+g(-x+ct) \right ]+
\frac{1}{2c}\,\left [\int_0^{x+ct} h(\xi)d\xi +\int_0^{ct-x} h(\xi)d\xi\right ] \qquad\nonumber \\ [6pt]
+\exp(\sigma (-x+ct)) \int_0^{x-ct} \exp(\sigma \xi)\,f(-\xi/c)\,d\xi \qquad \qquad \qquad \quad \qquad \nonumber \\ [6pt] - 2\sigma \exp(\sigma (-x+ct)) \int_0^{x-ct} \exp(\sigma \xi)\,\phi(-\xi)\,d\xi, ~ (x,t)\in D^{-}\,. \qquad
\end{eqnarray*}

Thus, the solution of the Robin problem (\ref{1-14}) in $D^{-}\cup D^{+}$ is
\begin{eqnarray}\label{1-19}
u(x,t)=
\left\{\begin{array}{ll}
\displaystyle \frac{1}{2}\,\left [ g(x+ct)+g(x-ct)\right ]	+\frac{1}{2c}\,\int_{x-ct}^{x+ct} h(\xi)d\xi,\,(x,t)\in D^{+},\\ [16pt]
\displaystyle \frac{1}{2}\,\left [g(x+ct)+g(-x+ct) \right ]+
\frac{1}{2c}\,\left [\int_0^{x+ct} h(\xi)d\xi +\int_0^{ct-x} h(\xi)d\xi\right ] \\ [16pt]
\displaystyle \qquad \qquad +\exp(\sigma (-x+ct)) \int_0^{x-ct} \exp(\sigma \xi)\,f(-\xi/c)\,d\xi   \\ [12pt]
\displaystyle \qquad \qquad  - 2\sigma \exp(\sigma (-x+ct)) \int_0^{x-ct} \exp(\sigma \xi)\,\phi(-\xi)\,d\xi, ~ (x,t)\in D^{-}\,.
\end{array}\right.
\end{eqnarray}

\section{Continuity properties of the solution across the characteristic line}

\subsection{Continuity properties of the Dirichlet problem solution}

From formula (\ref{1-10}) for the solution of the Dirichlet problem (\ref{1-4}) it follows that along the characteristic line $t=x/c$  this solution has the following form:
\begin{eqnarray}\label{1-22}
u(x,x/c)=
\left\{\begin{array}{ll}
\displaystyle \frac{1}{2}\,\left [ g(2x)+g(0)\right ]	+\frac{1}{2c}\,\int_{0}^{2x} h(\xi)d\xi,\,(x,t)\in D^{+},\\ [14pt]
\displaystyle  \frac{1}{2}\,\left [ g(2x)-g(0)\right ]	+\frac{1}{2c}\,\int_{0}^{2x} h(\xi)d\xi+\nu(0),\, (x,t)\in D^{-}\,.
\end{array}\right.
\end{eqnarray}
Hence,
\begin{eqnarray}\label{1-23}
[u]_{t=x/c}:=u(x,x/c)_{(x,t)\in D^{-}}- u(x,x/c)_{(x,t)\in D^{+}}=\nu(0)-g(0).
\end{eqnarray}

\noindent \textbf{Proposition 1.} \emph{The solution to the Dirichlet problem (\ref{1-4}) is discontinuous across the characteristic line $t=x/c$. The jump, as given by (\ref{1-23}), is determined by the initial data $g(0)$ and the boundary data $\nu(0)$ in (\ref{1-4}).}\\

\noindent \textbf{Corollary 1.} \emph{The solution to the Dirichlet problem (\ref{1-4}) is continuous across the characteristic line $t=x/c$ if the following compatibility condition is satisfied:}
\begin{eqnarray}\label{1-24}
g(0)=\nu(0).
\end{eqnarray}

In inverse problems where the corresponding direct problem includes a Dirichlet condition, the additional condition (measured output) is typically defined by a Neumann condition. Such inverse problems are associated with the Dirichlet-to-Neumann operator \cite{VGR-AH-3:2021}. Therefore, it is important to examine the continuity properties across the characteristic lines, not only for the solution itself but also for its first-order partial derivatives.\\

\noindent \textbf{Proposition 2.} \emph{The first partial derivatives $u_x$ and $u_t$ of the solution of the Dirichlet problem (\ref{1-4}) are discontinuous across the characteristic line $t=x/c$, with the jumps, given by
\begin{eqnarray}\label{1-25}
\displaystyle [u_x]_{t=x/c}:=u_x(x,x/c)_{(x,t)\in D^{-}}- u_x(x,x/c)_{(x,t)\in D^{+}}=\frac{1}{c}\,h(0)-\frac{1}{c}\,\nu'(0)
\end{eqnarray}
and
\begin{eqnarray}\label{1-26}
\displaystyle [u_t]_{t=x/c}:=u_t(x,x/c)_{(x,t)\in D^{+}}- u_t(x,x/c)_{(x,t)\in D^{-}}=-h(0)+\nu'(0)
\end{eqnarray}
respectively, across the characteristic being determined by the initial $h(0)$ and boundary $\nu(0)$ data in ((\ref{1-4})).} \\

\noindent \textbf{Corollary 2.} \emph{Under the compatibility condition
\begin{eqnarray*}
\displaystyle h(0)=\nu'(0)
\end{eqnarray*}
both first partial derivatives $u_x$ and $u_t$ of the solution of the Dirichlet problem (\ref{1-4}) are  continuous across the characteristic line $t=x/c$.}

\subsection{Continuity properties of the Neumann problem solution}

From formula (\ref{1-13}) it follows that along the characteristic line $t=x/c$ the solution of the Neumann problem (\ref{1-11}) has the following form:
\begin{eqnarray}\label{1-27}
u(x,x/c)=
\left\{\begin{array}{ll}
\displaystyle \frac{1}{2}\,\left [ g(2x)+g(0)\right ]	+\frac{1}{2c}\,\int_{0}^{2x} h(\xi)d\xi,\,(x,t)\in D^{+},\\ [14pt]
\displaystyle  \frac{1}{2}\,\left [ g(2x)+g(0)\right ]	+\frac{1}{2c}\,\int_{0}^{2x} h(\xi)d\xi,\, (x,t)\in D^{-}\,.
\end{array}\right.
\end{eqnarray}
Hence,
\begin{eqnarray}\label{1-28}
[u]_{t=x/c}=0.
\end{eqnarray}

\noindent \textbf{Proposition 3.} \emph{The solution to the Neumann problem (\ref{1-11}) is continuous across the characteristic line $t=x/c$, without requiring any compatibility condition on the initial and boundary data.}\\

\noindent \textbf{Proposition 4.} \emph{The first partial derivatives $u_x$ and $u_t$ of the solution of the Neumann problem (\ref{1-11}) are discontinuous across the characteristic line $t=x/c$, with the jumps, given by
\begin{eqnarray*}
\displaystyle [u_x]_{t=x/c}:=u_x(x,x/c)_{(x,t)\in D^{-}}- u_x(x,x/c)_{(x,t)\in D^{+}}=-g'(0)+f(0)
\end{eqnarray*}
and
\begin{eqnarray*}
\displaystyle [u_t]_{t=x/c}:=u_t(x,x/c)_{(x,t)\in D^{-}}- u_t(x,x/c)_{(x,t)\in D^{+}}=c g'(0)-c f(0)
\end{eqnarray*}
respectively, across the characteristic being determined by the initial $h(0)$ and boundary $\nu(0)$ data in ((\ref{1-4})).} \\

\noindent \textbf{Corollary 3.} \emph{Under the compatibility condition
\begin{eqnarray}\label{1-29}
g'(0)= f(0)
\end{eqnarray}
both first partial derivatives $u_x$ and $u_t$ of the solution of the Neumann problem (\ref{1-1}) are  continuous, respectively, across the characteristic line $t=x/c$.}

\subsection{Continuity properties of the Robin problem solution}

\noindent \textbf{Proposition 5.} \emph{The solution to the Robin problem (\ref{1-14}) is continuous across the characteristic line $t=x/c$, without requiring any compatibility condition on the initial and boundary data.}\\

\noindent \textbf{Proposition 6.} \emph{The first partial derivatives $u_x$ and $u_t$ of the solution of the Robin problem (\ref{1-14}) are discontinuous across the characteristic line $t=x/c$, with the jumps, given by
\begin{eqnarray*}
\displaystyle [u_x]_{t=x/c}:=u_x(x,x/c)_{(x,t)\in D^{-}}- u_x(x,x/c)_{(x,t)\in D^{+}}=-g'(0)-\sigma g(0)+f(0)
\end{eqnarray*}
and
\begin{eqnarray*}
\displaystyle [u_t]_{t=x/c}:=u_t(x,x/c)_{(x,t)\in D^{-}}- u_t(x,x/c)_{(x,t)\in D^{+}}=c g'(0)+\sigma c g(0)-c f(0)
\end{eqnarray*}
respectively.} \\

\noindent \textbf{Corollary 4.} \emph{Both first partial derivatives $u_x$ and $u_t$ of the solution of the Robin problem (\ref{1-14}) are continuous across the characteristic line $t=x/c$, provided that the compatibility condition
\begin{eqnarray*}
g'(0)+\sigma  g(0)= f(0)
\end{eqnarray*}
is satisfied.}

All of these assertions are derived directly from the solution of the Robin problem (\ref{1-14}), defined by formula (\ref{1-19}).

\end{document}